\documentclass[a4paper,10pt]{article} 
\usepackage{amsmath} 
\usepackage{graphicx} 
\usepackage{amsfonts} 
\usepackage{amssymb} 
\usepackage[latin1]{inputenc}

\newcommand{\la}{\lambda} 

\newcommand{\cc}{{\cal C}}

\newcommand{\cf}{{\cal F}}

\newcommand{\cl}{{\cal L}}

\newcommand{\si}{\sigma}

\newcommand{\R}{{\mathbb R}} 
 
\newcommand{\G}{{\mathbb P}}

\newtheorem{theorem}{Theorem}[section]

\newtheorem{definition}[theorem]{Definition}

\newtheorem{lemma}[theorem]{Lemma}

\newtheorem{proposition}[theorem]{Proposition}

\begin{document}
\thispagestyle{empty}

\begin{center}
\huge {

Delay equations with non-negativity constraints 
driven by a H\"older continuous function of order $\beta \in (\frac13,\frac12)$}

\vspace{.5cm}

\normalsize {\bf Mireia Besal\'u}, {\bf David M\'arquez-Carreras$^*$} and {\bf Carles
Rovira}

{\footnotesize \it

Facultat de Matem\`atiques, Universitat de Barcelona, Gran
Via 585, 08007-Barcelona.

 {\it E-mail addresses}: mbesalu@ub.edu, davidmarquez@ub.edu, 
carles.rovira@ub.edu}

{$^*$ Corresponding author.}

\end{center}

\begin{abstract}%
In this note we prove an existence and uniqueness result of
solution for multidimensional  delay differential  equations with normal reflection and driven by a H\"older continuous function of order $\beta \in (\frac13,\frac12)$. We also obtain a bound for the supremum norm of this solution. As an application, we get  these results for stochastic differential equations  driven by a fractional Brownian motion with Hurst
parameter H $\in (\frac13,\frac12)$.

\end{abstract}

\vspace{0.1cm}


\noindent {\bf AMS 2000 MSC:} 60H05, 60H07

\noindent {\bf Running head:}  delay equations with constraints

\setcounter{section}{0}

%%%%%%%%%%%%%%%%%%%%%%%%%%%%%%%%%%%%%%%%%%%%%%%%%%%%%%%%%%%%%%%%%%%%%%%%%%%%

\renewcommand{\theequation}{1.\arabic{equation}}
\setcounter{equation}{0}
\section{Introduction}

The theory of rough path analysis has been developed from the initial paper
by Lyons \ \cite{L}. The aim of this theory is to analyze dynamical
systems $dx_{t}=f(x_{t})dy_{t}$, where the control function $y$ is not
differentiable but has finite $p$-variation for some $p> 1$.
There is a wide literature on rough path analysis (see, for instance,
Lyons and Qian \cite{L-Q},  Friz and Victoir  \cite{FV2}, Lejay \cite{Le}, Lyons   \cite{LCL} or  Gubinelli  \cite{Gu}).

 A path-wise approach to classical stochastic
calculus has been one of the motivations to build rough path analysis theory. A nice application of the rough path
analysis is the stochastic calculus with respect to the fractional Brownian
motion with Hurst parameter $H\in (0,1)$. We
refer, for instance to Coutin and Lejay \cite{CA}, Friz and Victoir \cite{Fr}, Friz \cite{FV} and Ledoux {\it et al.} \cite{LQZ} for some
applications of   rough path analysis to the stochastic calculus.

 Nualart and R\u{a}\c{s}canu in  \cite{NR} developed an alternative approach to the study of dynamical
systems   $dx_{t}=f(x_{t})dy_{t}$, where the control function $y$ is H\"{o}lder
continuous of order $\beta >\frac{1}{2}$.
In this case the Riemann-Stieltjes integral $\int_{0}^{t}f(x_s)dy_s$     can be expressed as a Lebesgue integral
using fractional derivatives  following the ideas of Z\"{a}hle \cite{Z}.
Later, Hu and Nualart \cite{H-N} extended this approach  to the case  $\beta \in (\frac{1}{3},\frac{1}{2})$ . In this work they give an  explicit expression for  the integral $\int_{0}^{t}f(x_s)dy_s$
that depends on the functions $x$, $y$ and a quadratic multiplicative functional 
$x\otimes y$.  Using this formula,  the authors have established the existence and uniqueness of a solution for the dynamical system $dx_{t}=f(x_{t})dy_{t}$ driven by a H\"{o}lder
continuous function $y$ of order $\beta \in (\frac{1}{3},\frac{1}{2})$. Finally,  using the same approach, Besal\'u and Nualart \cite{B-N} got estimates
for the supremum norm of the solution.

 The purpose of this paper is to study a differential delay equation with non-negativity constraints 
 driven by a H\"older continuous function $y$ of order $%
\beta\in\left(\frac{1}{3},\frac{1}{2}\right)$ using the methodology introduced in \cite{H-N}.
 We will consider the problems of existence, uniqueness and boundedness of the solutions.  As an
application we will study a
stochastic delay  differential equations with non-negativity constraints  driven by a fractional Brownian motion
with Hurst parameter $H\in\left(\frac{1}{3},\frac{1}{2}\right)$. 
These results extend the work by Besal\'u and Rovira \cite{B-R}, where is considered the case $H>\frac 12$.

More precisely, we consider a  delay differential equation with positivity constraints on $\R^d$ of the form:
\begin{eqnarray}
x(t)&=&\eta(0)+\int_0^t b(s,x)ds+\int_0^t \si(x(s-r))dy_s+z(t),\quad t\in(0,T],\nonumber\\
x(t)&=& \eta(t),\qquad t\in[-r,0], \nonumber
\end{eqnarray} 
where $r$ denotes a strictly positive time delay, $y$ is a $m$-dimensional $\beta$-H\"older continuous function  with  $\frac13<\beta<\frac{1}{2}$, $b(s,x)$ the hereditary term, depends on the path $\left\{x(u),-r\leq u\leq s\right\}$, while $\eta:[-r,0]\rightarrow\R^d_+$ is a non negative smooth function, with $\R^d_+=\left\{u\in\R^d;\,u_i\geq 0\;\mathrm{for}\;i=1,\ldots,d\right\}$ and $z$ is a vector-valued non-decreasing process which ensures that the non-negativity constraints on $x$ are enforced.

\vskip 7pt

Then, we  will apply pathwise  our deterministic result to  a stochastic delay differential equation with positivity constraints on $\R^d$ of the form:
\begin{eqnarray}
X(t)&=&\eta(0)+\int_0^t b(s,X)ds+\int_0^t \si(X(s-r))dW_s^H+Z(t),\quad t\in(0,T],\nonumber\\
X(t)&=& \eta(t),\qquad t\in[-r,0], \nonumber
\end{eqnarray} 
where  $W^H=\left\{W^{H,j},\,j=1,\ldots,m\right\}$ are independent fractional Brownian motions with Hurst parameter $\frac13<H<\frac{1}{2}$ defined in a complete probability space $(\Omega,\cf,\G)$,  while $\eta:[-r,0]\rightarrow\R^d_+$ is a deterministic  non negative smooth function and $Z$ is a vector-valued non-decreasing process which ensures that the non-negativity constraints on X are enforced.

\vskip 4pt

As far as we know, stochastic delay differential equations with constraints and driven by a fractional Brownian motion has only been considered when $H>\frac12$ (\cite{B-R}). Furthermore, 
the literature about stochastic delay differential equations  driven by a fractional Brownian
motion is scarce. For the  case $H > \frac12$ has been studied the existence and uniqueness of solution (\cite{FR}, \cite{LT}),
the existence and regularity of the density (\cite{LT}) and the  convergence when the delay goes to zero (\cite{F-R}). For $H<\frac12$ we can find the
results about the existence and uniqueness of solution  (\cite{NNT}, \cite{TT}). Actually, in \cite{NNT}  the authors consider a similar equation to our case  but without reflection. Moreover, they use another approach in order to define the stochastic integral based on L\'evy area. In any case,
we will use some results on fractional Brownian motion taken from this paper.

\vskip 4pt

Anyway, as it has been described in this paper of Kinnally and Williams \cite{K-W} there are some models afected by some type of noise where the dynamics are related to propagation delay and some of them are naturally non-negative quantities. So, it is natural to continue the study of the  stochastic delay differential equations and non-negativity constraints driven by a fractional Brownian motion.

\vskip 7pt

In our work, we will make use of the techniques introduced by Hu and Nualart \cite{H-N}  with some ideas borrowed from Besal\'u and Rovira \cite{B-R}. In this framework, let us point out again that one novelty of our paper is the non-negative constraints dealing with equations driven by a H\"older continuous function of order $\beta \in (\frac13,\frac12)$. We have used the
Skorohod's  mapping. Let us recall now the Skorokhod problem.
Set  $$\cc_+(\R_+,\R^d):=\left\{x\in\cc(\R_+,\R^d): x(0)\in\R^d_+\right\}.$$
\vskip 5pt
\noindent

\begin{definition}
Given a path $z\in\cc_+(\R_+,\R^d)$, we say that a pair $(x,y)$ of functions in $\cc_+(\R_+,\R^d)$ solves the Skorokhod problem for $z$ with reflection if
\begin{enumerate}
\item $x(t)=z(t)+y(t)$ for all $t\geq 0$ and $x(t)\in\R_+^d$ for each $t\geq 0$,
\item for each $i=1,\ldots, d$, $y^i(0)=0$ and $y^i$ is nondecreasing,
\item for each $i=1,\ldots, d$, ${\displaystyle \int_0^t x^i(s)dy^i_s=0}$ for all $t\geq 0$, so $y^i$ can increase only when $x^i$ is at zero.
\end{enumerate}
\end{definition}
\noindent
It is known that we have an explicit formula for $y$ in terms of $z$: for each $i=1,\ldots,d$
\begin{equation*}
y^i(t)=\max_{s\in[0,t]} \left(z^i(s)\right)^-.
\end{equation*}

The path $z$ is called the reflector of $x$ and the path $y$ is called the regulator of $x$. We use the Skorokhod mapping for constraining a continuous real-valued function to be non-negative by means of reflection at the origin. 

\vskip 7pt

The structure of the paper is as follows:
in the next section we give some preliminaries, our hypothesis and we state the main results of our paper.
In Section 3, we give some basic facts about fractionals integrals. 
Section 4  is devoted to prove our main result: the existence and
uniqueness  for the solution for deterministic equations, while Section 5 deals with the problem of the boundedness.
In  Section 6 we  apply the deterministic results to the stochastic case. Finally, Section 7
is devoted to give some technical results, as a fixed point theorem, and some properties related to the Skorohod problem.

%%%%%%%%%%%%%%%%%%%%%%%%%%%%%%%%%%%%%%%%%%%%%%%%%%%%%%%%%%%%%%%%%%%
\renewcommand{\theequation}{2.\arabic{equation}}
\setcounter{equation}{0}
\section{Main results}
\vskip 5pt
\noindent
Fix a time interval $[0,T]$. For any function $x:[0,T]\rightarrow \mathbb{R}%
{}^{n}$, the $\gamma $-H\"{o}lder norm of $x$ on the interval $[s,t] \subset \lbrack 0,T]$, where $0< \gamma \le 1$, will be denoted by 
\begin{equation*}
\left\| x\right\| _{\gamma(s,t) }=\sup_{s< u<v< t}\frac{|x(v)-x(u)|}{%
(v-u)^{\gamma }}.
\end{equation*}%
If $\Delta_T :=\left\{ (s,t):0\leq s<  t\leq T\right\} $, for any $(s,t)\in
\Delta_T $ and for any $g:\Delta_T \rightarrow {\mathbb{R}}^{n}$ we set 
\begin{equation*}
\left\| g\right\| _{\gamma(s,t)}=\sup_{s< u<v< t}\frac{|g(u,v)|}{%
(v-u)^{\gamma }}.
\end{equation*}%
We will also set $\left\| x\right\| _{\gamma}=\left\| x\right\|
_{\gamma(0,T)}$ and $\left\| x\right\| _{\gamma(r)}=\left\| x\right\|
_{\gamma(-r,T)}$. Moreover, $\left\| \cdot \right\| _{\infty(s,t)}$ will
denote the supremum norm in the interval $(s,t)$, and for simplicity $\left\| x\right\| _{\infty}=\left\| x\right\|
_{\infty(0,T)}$ and $\left\| x\right\| _{\infty(r)}=\left\| x\right\|
_{\infty(-r,T)}$.
\vskip 5pt
\noindent
Fix $0<\beta \le 1$. As in \cite{L} we introduce the
following definition.

\begin{definition}\label{d2.1}
We will say that $(x,y,x\otimes y)$ is an {$(d,m)$-dimensional $\beta $-H%
\"{o}lder continuous multiplicative functional} if:

\begin{enumerate}
\item $x:[0,T]\rightarrow \mathbb{R}{}^{d}$ and $y:[0,T]\rightarrow \mathbb{R%
}{}^{m}$ are $\beta $-H\"{o}lder continuous functions,

\item $x\otimes y:\Delta_T \rightarrow \mathbb{R}{}^{d}\otimes \mathbb{R}{}^{m}
$ is a continuous function satisfying the following properties:

\begin{enumerate}
\item {\upshape(Multiplicative property)} For all $s\leq u\leq t$ we have 
\begin{equation*}
(x\otimes y)_{s,u}+(x\otimes y)_{u,t}+(x(u)-x(s))\otimes
(y(t)-y(u))=(x\otimes y)_{s,t}.
\end{equation*}

\item For all $(s,t)\in \Delta_T $ 
\begin{equation*}
\left| (x\otimes y)_{s,t}\right| \leq c|t-s|^{2\beta }.
\end{equation*}
\end{enumerate}
\end{enumerate}
\end{definition}
We will denote by $M_{d,m}^{\beta
}(0,T)$ the space of $(d,m)$-dimensional $\beta $-H\"{o}lder continuous
multiplicative functionals. 
Furthermore, we will denote by $M_{d,m}^{\beta
}(a,b)$ the obvious extension of the definition $M_{d,m}^{\beta
}(0,T)$ to a general interval $(a,b)$. 
We refer the reader to \cite{L} and  \cite{H-N} for a more detailed presentation on $\beta $-H\"{o}lder continuous
multiplicative functionals.
\vskip 5pt

 We have now the tools to give our results.  Set $\beta\in(\frac13,\frac12)$. 

Consider   the deterministic stochastic differential equation on $\R^d$
\begin{eqnarray}
x(t)&=&\eta(0)+\int_0^t b(s,x)ds+\int_0^t\sigma(x(s-r))dy_s+ z(t),\quad t\in (0,T],\nonumber\\
x(t) &=& \eta(t),\quad t\in[-r,0],\label{det}
\end{eqnarray} 
{where for each} $i=1,\ldots,d$,
\begin{equation*}
{\mathit z^i(t)=\max_{s\in[0,t]}\left(\xi^i(s)\right)^-,\qquad t\in[0,T],}
\end{equation*}
{and}
\begin{equation*}
{\mathit \xi(t)= \eta(0)+\int_0^t b(s,x)ds+\int_0^t \si(x(s-r))dy_s,\qquad t\in[0,T].}
\end{equation*}
Let us consider the following hypothesis:
\begin{itemize}
\item[\bfseries(H1)] $\sigma:\R^d\rightarrow \R^d\times\R^m$ is a continuously differentiable function such that $\sigma'$ is locally $\gamma$- H\"older continuous for $\gamma>\frac{1}{\beta}-2$.
\item[\bfseries(H2)] $b:[0,T]\times C(-r,T;\R^d)\rightarrow \R^d$ is a measurable function such that for every $t>0$ and $f\in C(-r,T;\R^d),\, b(t,f)$ depends only on $\left\{f(s);-r\leq s\leq t\right\}$. Moreover, there exists $b_0\in L^\rho(0,t;\R^d)$ with $\rho\geq 2$ and $\forall N\geq 0$ there exists $L_N>0$ such that:
\begin{itemize}
\item[{\bf (1)}] $\left|b(t,x)-b(t,y)\right|\leq L_N \left\|x-y\right\|_{\infty(-r,t)},\; \forall x,y\;\mathrm{such\,that} \left\|x\right\|_{\infty(r)}\leq N,$ 
$\left\|y\right\|_{\infty(r)}\leq N,\;\forall t\in[0,T]$,
\item[{\bf (2)}] $\left|b(t,x)\right|\leq L_0 \left\|x\right\|_{\infty(-r,t)}+b_0(t),\quad \forall t\in[0,T]$.
\end{itemize}
\end{itemize}

\vskip 5pt
\noindent
The result of existence and uniqueness states as follows:

\begin{theorem}\label{tt1}
Assume that $\sigma$ and $b$ satisfy the hypothesis {\upshape \bfseries(H1)} and  {\upshape \bfseries(H2)} respectively with $\rho\geq \frac{1}{1-\beta}$. Assume also  that $\eta\geq 0$, $(\eta_{\cdot-r},y,\eta_{\cdot-r}\otimes y)\in M_{d,m}^{\beta}(0,r)$ and $(y_{\cdot-r},y,y_{\cdot-r}\otimes y)\in M^\beta_{m,m}(r,T)$. Then the equation (\ref{det}) has a unique solution $x\in \cc(-r,T;\R^d_+)$.
\end{theorem}

\vskip 5pt

\noindent {\bf Remark}. If we assume that $\eta\geq 0$ is a differentiable continuous function with positive derivative, then the assumptions on $\eta$ of this theorem are satisfied.

\medskip

\noindent
In order to study the boundedness of the solutions we need to stronger our hyphotesis. Consider now:
\begin{itemize}
%\item[\bfseries(H1)] $\sigma:\R^d\rightarrow \R^d\times\R^m$ is a continuously differentiable %function such that $\sigma'$ is bounded and  locally $\gamma$- H\"older continuous for %$\gamma>\frac{1}{\beta}-2$
\item[\bfseries(H3)] $b$ and $\sigma'$ are bounded function.
\end{itemize}
Then, the result is as follows:
\begin{theorem}\label{tt2}
Assume that $\sigma$ and $b$ satisfy the hypothesis {\upshape \bfseries(H1)}, {\upshape \bfseries(H2)} and  {\upshape \bfseries(H3)}. Also assume that $\eta\geq 0$  satisfies $(\eta_{\cdot-r},y,\eta_{\cdot-r}\otimes y)\in M_{d,m}^{\beta}(0,r)$ and finally that $(y_{\cdot-r},y,y_{\cdot-r}\otimes y)\in M^\beta_{m,m}(r,T)$. Set
$$\mu=\left\|b\right\|_{\infty}+\left\|\sigma\right\|_{\infty}+\left\|\sigma'\right\|_{\infty}+\left\|\sigma'\right\|_{\gamma}.$$
Then, the solution of (\ref{det}) is bounded as follows
\begin{equation}
\|x\|_{\infty}\le 2+ \eta(0) +  T \left\{K\left(\left\|\eta\right\|_{\beta}
+\left\|\eta_{\cdot-r}\otimes y\right\|_{2\beta}
+\mu (d^\frac12 +1)  \left[\left\|y\right\|_{\beta}
+\left\|y\right\|_{\beta}^2+ \left\|y_{\cdot-r}\otimes y\right\|_{2\beta}\right]\right)\right\}^\frac{1}{\beta},\label{e5.1}\end{equation}
where $K$ is a universal constant depending only on $\beta$ and  $\gamma$,
and
\begin{eqnarray*}
\left\|\eta\right\|_{\beta}&:=&\left\|\eta\right\|_{\beta(-r,0)},\\
\left\|\eta_{\cdot-r}\otimes y\right\|_{2\beta}
&:=&\left\|\eta_{\cdot-r}\otimes y\right\|_{2\beta(0,r)},\\
\left\|y\right\|_{\beta}&:=&\left\|y\right\|_{\beta(0,T)},\\
\left\|y_{\cdot-r}\otimes y\right\|_{2\beta}
&:=&
\left\|y_{\cdot-r}\otimes y\right\|_{2\beta(r,T)}.
\end{eqnarray*}
\end{theorem}

\vskip 5pt
\noindent
Our last result is an application of the above theorems to stochastic delay differential equations. More precisely, let us consider a
stochastic delay differential equation with positivity constraints on $\R^d$ of the form:
\begin{eqnarray}
X(t)&=&\eta(0)+\int_0^t b(s,X)ds+\int_0^t \si(X(s-r))dW_s^H+Z(t),\quad t\in(0,T],\nonumber\\
X(t)&=& \eta(t),\qquad t\in[-r,0], \label{eqstoc}
\end{eqnarray} 
where  $W^H=\left\{W^{H,j},\,j=1,\ldots,m\right\}$ are independent fractional Brownian motions with Hurst parameter $\frac13<H<\frac{1}{2}$ defined in a complete probability space $(\Omega,\cf,\G)$, 
 {and for each} $i=1,\ldots,d$
\begin{equation*}
{\mathit Z^i(t)=\max_{s\in[0,t]}\left(\Xi^i(s)\right)^-,\qquad t\in[0,T],}
\end{equation*}
{and}
\begin{equation*}
{\mathit \Xi(t)= \eta(0)+\int_0^t b(s,X)ds+\int_0^t \si(s,X(s-r))dW^H_s,\qquad t\in[0,T].}
\end{equation*}

\noindent
Then, our result for the stochastic equation, reads as follows:

\begin{theorem}\label{tt3} Fix $\beta \in (\frac13,H)$.
Assume that $\sigma$ and $b$ satisfy the hypothesis {\upshape \bfseries(H1)} and  {\upshape \bfseries(H2)} respectively with $\rho\geq \frac{1}{1-\beta}$. Assume also  that $\eta$ is a non-negative bounded function such that $(\eta_{\cdot-r},W^H,\eta_{\cdot-r}\otimes W^H)\in M_{d,m}^{\beta}(0,r)$ almost surely. Then the equation (\ref{eqstoc}) has a unique solution
 \[X\in L^0(\Omega,\cf,\G; \cc(-r,T;\R^d_+)).\]

\noindent
Furthermore, if  {\upshape \bfseries(H3)} is satisfied and $E(\left\|\eta_{\cdot-r}\otimes W^H\right\|_{2\beta(0,r)}^p)  < \infty, 
\forall p \ge 1$, then $E ( \Vert X \Vert^p_{\infty} ) < \infty,
\forall p \ge 1.$
\end{theorem}

%%%%%%%%%%%%%%%%%%%%%%%%%%%%%%%%%%%%%%%%%%%%%%%%%%%%%%%%%%%%%%%%%%%%%%%%%%%
\renewcommand{\theequation}{3.\arabic{equation}} \setcounter{equation}{0}

\section{Fractional integrals and derivatives}

In this section we recall some definitions and results on fractional integrals. We refer the reader to \cite{H-N} for a more detailed presentation.
\vskip 5pt
\noindent
Let $a,b\in  \mathbb{R}$ with $a<b$. Let $f\in L^{1}(a,b)$ and $\alpha >0$. The
left-sided and right-sided fractional Riemann-Liouville integrals of $f$ of
order $\alpha $ are defined for almost all $t\in (a,b)$ by 
\begin{equation*}
I_{a+}^{\alpha }f(t)=\frac{1}{\Gamma (\alpha )}\int_{a}^{t}(t-s)^{\alpha
-1}f(s)ds
\end{equation*}%
and 
\begin{equation*}
I_{b-}^{\alpha }f(t)=\frac{(-1)^{-\alpha }}{\Gamma (\alpha )}%
\int_{t}^{b}(s-t)^{\alpha -1}f(s)ds,
\end{equation*}%
respectively, where $(-1)^{-\alpha }=e^{-i\pi \alpha }$ and $\Gamma (\alpha
)=\int_{0}^{\infty }r^{\alpha -1}e^{-r}dr$ is the Euler gamma function. For any $p\ge 1$, let $%
I_{a+}^{\alpha }(L^{p})$ (resp. $I_{b-}^{\alpha }(L^{p})$) be the image of $%
L^{p}(a,b)$ by the operator $I_{a+}^{\alpha }$ (resp. $I_{b^{-}}^{\alpha }$%
). If $f\in I_{a+}^{\alpha }(L^{p})$ (resp. $f\in I_{b-}^{\alpha }(L^{p})$)
and $0<\alpha <1$, then the Weyl derivatives are defined as 
\begin{eqnarray}
D_{a+}^{\alpha }f(t) &=&\frac{1}{\Gamma (1-\alpha )}\left( \frac{f(t)}{%
(t-a)^{\alpha }}+\alpha \int_{a}^{t}\frac{f(t)-f(s)}{(t-s)^{\alpha +1}}%
ds\right) ,  \nonumber \\
D_{b-}^{\alpha }f(t) &=&\frac{(-1)^{\alpha }}{\Gamma (1-\alpha )}\left( 
\frac{f(t)}{(b-t)^{\alpha }}+\alpha \int_{t}^{b}\frac{f(t)-f(s)}{%
(s-t)^{\alpha +1}}ds\right) , \nonumber% \label{d3}
\end{eqnarray}%
where $a\leq t\leq b$ (the convergence of the integrals at the singularity $%
s=t$ holds point-wise for almost all $t\in (a,b)$ if $p=1$ and moreover in
the $L^{p}$-sense if $1<p<\infty $).\newline
\noindent If $f\in C^{\lambda }(a,b)$ and $g\in C^{\mu }(a,b)$ with $\lambda
+\mu >1$, it is proved in \cite{Z} that the Riemman-Stieltjes integral $%
\int_{a}^{b}fdg$ exists. The following proposition provides an explicit
expression for the integral $\int_{a}^{b}fdg$ in terms of fractional
derivatives (see \cite{Z}).

\begin{proposition}
Suppose that $f\in C^{\lambda }(a,b)$ and $g\in C^{\mu }(a,b)$ with $\lambda
+\mu >1$. Let $1-\mu <\alpha <\lambda $. Then the Riemann-Stieltjes integral 
$\int_{a}^{b}fdg$ exists and it can be expressed as 
\begin{equation}
\int_{a}^{b}fdg=(-1)^{\alpha }\int_{a}^{b}D_{a+}^{\alpha
}f(t)D_{b-}^{1-\alpha }g_{b-}(t)dt,  \label{forpart}
\end{equation}%
where $g_{b-}(t)=g(t)-g(b).$
\end{proposition}
\noindent
But if $x,y \in C^{\beta} (a,b)$ with $\beta \in (\frac13,\frac12)$ we can not use Equation (\ref{forpart}) to define the integral $\int_{a}^{b}f(x(t))dy_{t}$, so
we need to recall the construction of the integral $\int_{a}^{b}f(x(t))dy_{t}$
given by Hu and Nualart in \cite{H-N} using fractional derivatives.
\begin{definition}
\label{defhn} Let $(x,y,x\otimes y)\in M_{d,m}^{\beta }(0,T)$. Let $f:%
\mathbb{R}{}^{d}\rightarrow \mathbb{R}{}^{m}\otimes \mathbb{R}{}^{d}$ be a
continuously differentiable function such that $f^{\prime }$ is locally $%
\lambda $-H\"{o}lder continuous, where $\lambda >\frac{1}{\beta }-2$. Fix $%
\alpha >0$ such that $1-\beta <\alpha <2\beta $, and $\alpha <\frac{\lambda
\beta +1}{2}$. Then, for any $0\leq a<b\leq T$ we define 
\begin{eqnarray*}
\int_{a}^{b}f(x(r))dy_{r}&=&(-1)^{\alpha
}\sum_{j=1}^{m}\int_{a}^{b}\widehat{D}_{a+}^{\alpha }f_{j}(x)(r)D_{b-}^{1-\alpha
}y_{b-}^{j}(r)dr \\
&& -(-1)^{2\alpha
-1}\sum_{i=1}^{d}\sum_{j=1}^{m}\int_{a}^{b}D_{a+}^{2\alpha -1}\partial
_{i}f_{j}(x)(r)D_{b-}^{1-\alpha }\mathcal{D}_{b-}^{1-\alpha }(x\otimes
y)^{i,j}(r)dr,
\end{eqnarray*}
\end{definition}
where for $r\in \left( a,b\right)$ 
\begin{equation*}
\begin{array}{l}
\displaystyle\widehat{D}_{a+}^{\alpha }f(x)(r)=\frac{1}{\Gamma (1-\alpha )}%
\left( \frac{f(x(r))}{(r-a)^{\alpha }}
%\right.  \\[4mm]\displaystyle\qquad \qquad \left.
 +\alpha \int_{a}^{r}\frac{%
f(x(r))-f(x(\theta ))-\sum_{i=1}^{m}\partial _{i}f(x(\theta))(x^{i}(r)-x^{i}(\theta))}{(r-\theta )^{\alpha +1}}d\theta \right) %
\end{array}%
\end{equation*}
is  the \textit{compensated fractional
derivative}  
and  
\begin{equation*}
\mathcal{D}_{b-}^{1-\alpha }(x\otimes y)(r)=\frac{(-1)^{1-\alpha }}{\Gamma
(\alpha )}\left( \frac{(x\otimes y)_{r,b}}{(b-r)^{1-\alpha }}+(1-\alpha
)\int_{r}^{b}\frac{(x\otimes y)_{r,s}}{(s-r)^{2-\alpha }}ds\right)
\end{equation*}
is the extension of the fractional derivative of $x\otimes y$.

\vskip 5pt
\noindent
Let us finish this section recalling two propositions from \cite{H-N}. In the sequel, $k$ denotes a generic constant that may depend on the
parameters $\beta, \alpha$ and $\gamma$.

\begin{proposition}\label{prop23}
Let $(x,y,x \otimes y)$ be in $M_{d,m}^\beta (0,T)$. Assume that $f: \R^d \longrightarrow \R^m$ is a continuous
differentiable function such that $f'$ is bounded and $\gamma$-H\"older continuous, where $\gamma > \frac{1}{\beta} -2.$ Then for any
$0 \le a < b \le T$ we have
\begin{eqnarray*}
\Vert \int f(x(r)) dy_r \Vert_{\beta(a,b)} &\le & k \vert f(x(a)) \vert \Vert y \Vert_{\beta(a,b)} + k \Phi_{a,b,\beta} (x,y)
\\
&&  \times \left( \Vert f' \Vert_\infty + \Vert f' \Vert_\gamma \Vert x \Vert^\gamma_{\beta(a,b)} (b-a)^{\gamma \beta} \right) (b-a)^\beta,
\end{eqnarray*}
where
$$ \Phi_{a,b,\beta} (x,y)= \Vert x \otimes y \Vert_{2\beta(a,b)} + \Vert x \Vert_{\beta(a,b)} \Vert y \Vert_{\beta(a,b)}.$$
\end{proposition}

\begin{proposition}\label{prop24}
Suppose that  $(x,y,x \otimes y)$ and   $(y,z,y \otimes z)$   belong to  $M_{d,m}^\beta (0,T)$. Let $f: \R^d \longrightarrow \R^m$ be a continuously
differentiable function such that $f'$ is $\gamma$-H\"older continuous and bounded, where $\gamma > \frac{1}{\beta} -2.$ 
Fix $\alpha >0$ such that $1 - \beta < \alpha <2 \beta, \alpha < \frac{\gamma \lambda +1}{2}.$ Then the following estimate holds:
\begin{eqnarray*}
\vert \int_a^b f(x(r)) d (y\otimes z)_{\cdot,b}(r) \vert & \le & k \vert f(x(a)) \vert\Phi_{a,b,\beta}(y,z) (b-a)^{2\beta}
\\
&& + k \left( \Vert f' \Vert_\infty + \Vert f' \Vert_\gamma \Vert x \Vert^\gamma_{\beta(a,b)} (b-a)^{\gamma \beta} \right) \Phi_{a,b,\beta}(x,y,z) (b-a)^{3\beta},
\end{eqnarray*}
where
$$ \Phi_{a,b,\beta} (x,y,z)= \Vert x \Vert_{\beta(a,b)} \Vert y \Vert_{\beta(a,b)} \Vert z \Vert_{\beta(a,b)}  + \Vert z \Vert_{\beta(a,b)}\Vert x \otimes y \Vert_{2\beta(a,b)} + \Vert x \Vert_{\beta(a,b)} \Vert y \otimes z  \Vert_{2\beta(a,b)}.$$
\end{proposition}

%%%%%%%%%%%%%%%%%%%%%%%%%%%%%%%%%%%%%%%%%%%%%%%%%%%%%%%%%%%%%%%%%%%%%%%%%%%%%%%%%%%%%%%%%%%%%%%
\renewcommand{\theequation}{4.\arabic{equation}}
\setcounter{equation}{0}
\section{Existence and uniqueness for deterministic integral equations}
\vskip 5pt
\noindent
The aim of this section is the proof of Theorem \ref{tt1}. For simplicity let us assume $T=Mr$.
\vskip 5pt
\noindent
{\bf Proof of Theorem \ref{tt1}:}
In order to prove that equation (\ref{det}) admits a unique continuous solution on $[-r,T]$, we will use an induction argument.
We shall prove that if the equation (\ref{det}) admits a unique solution $x^{(n)}$  on $[-r,nr]$ we can  prove that there is a unique solution $x^{(n+1)}$ on $[-r,(n+1)r]$. More precisely, 
our induction hypothesis is the following:\\
$\bf{(H_n)}$ The equation
\begin{eqnarray*}
x^{(n)}(t)&=&\eta(0)+\int_0^t b(s,x^{(n)})ds+\int_0^t\sigma(x^{(n-1)}(s-r))dy_s+ z^{(n)}(t),\quad t\in (0,nr],\nonumber\\
x^{(n)}(t) &=& \eta(t),\quad t\in[-r,0],\label{sn}
\end{eqnarray*}
where for $i=1,\dots,d$, $(z^{(n)})^{i}(t)=\max_{s\in[0,t]}((\xi^{(n)})^{i}(s))^-$ with
\[\xi^{(n)}(t)=\eta(0)+\int_0^t b(s,x^{(n)})ds+\int_0^t\sigma(x^{(n-1)}(s-r))dy_s,\]
has a unique solution $x^{(n)}\in \cc(-r,nr,\R_+^d)$ and moreover $(x^{(n)}_{\cdot-r},y,x^{(n)}_{\cdot-r}\otimes y)\in M^\beta_{d,m}(0,(n+1)r)$.
\vskip 5pt
\noindent
Actually, when we want to check  $\bf{(H_{n+1})}$ assuming $\bf{(H_{n})}$, we can write the equation of $\bf{(H_{n+1})}$
 as
\begin{eqnarray}
x^{(n+1)}(t)&=&\eta(0)+\int_0^t b(s,x^{(n+1)})ds+\int_0^t\sigma(x^{(n)}(s-r))dy_s+ z^{(n+1)}(t),\quad t\in (0,(n+1)r],\nonumber\\
x^{(n+1)}(t) &=& \eta(t),\quad t\in[-r,0].\label{sn1}
\end{eqnarray}
Since $(x^{(n)}_{\cdot-r},y,x^{(n)}_{\cdot-r}\otimes y)\in M^\beta_{d,m}(0,(n+1)r)$ we know that we can use  Definition {\ref{defhn}}  to define the integral $\int_0^t\sigma(x^{(n)}(s-r))dy_s$ appearing in  equation (\ref{sn1}).
Then, 
the proof will consist in checking the following  steps:
\begin{enumerate}
\item Existence of a solution of the equation (\ref{sn1}) in the space $\cc(-r,(n+1)r;\R_+^d)$.
\item Uniqueness of a solution of the equation (\ref{sn1}) in the space $\cc(-r,(n+1)r;\R_+^d)$.
\item The solution $x^{(n+1)}$ satifies that  $(x^{n+1}_{\cdot-r},y,x^{(n+1)}_{\cdot-r}\otimes y)\in M_{d,m}^{\beta}(0,(n+2)r)$.
\end{enumerate}

Actually, we will only proof the first case, that is  $\bf{(H_{1})}$. Notice that 
the induction step, that is the proof of $\bf{(H_{n+1})}$ assuming that $\bf{(H_{n})}$ is true, can be done repeating  the computations of this initial case.
\vskip 5pt
\noindent So, let us check   $\bf{(H_{1})}$. We will deal with the equation
\begin{eqnarray}
x^{(1)}(t)&=&\eta(0)+\int_0^t b(s,x^{(1)})ds+\int_0^t\sigma(\eta(s-r))dy_s+ z^{(1)}(t),\quad t\in (0,r],\nonumber\\
x^{(1)}(t) &=& \eta(t),\quad t\in[-r,0],\label{s1}
\end{eqnarray}
where for $i=1,\;\ldots,\;d$,  $(z^{(1)})^i(t)=\max_{s\in[0,t]}((\xi^{(1)})^i(s))^-$ and
\[\xi^{(1)}(t)=\eta(0)+\int_0^t b(s,x^{(1)})ds+\int_0^t\sigma(\eta(s-r))dy_s.\]
Note that since $(\eta_{\cdot-r},y,\eta_{\cdot-r}\otimes y)\in M_{d,m}^{\beta}(0,r)$ we can use Definition \ref{defhn} in order to define the  integral $\int_0^t\sigma(\eta(s-r))dy_s$ appearing in (\ref{s1}). The proof of this initial case will be divided en 3 steps:
\begin{enumerate}
\item Existence of a solution in the space $\cc(-r,r;\R^d_+)$.
\item Uniqueness of a solution in the space $\cc(-r,r;\R^d_+)$.
\item The solution $x^{(1)}$ satisfies that $(x^{(1)}_{\cdot-r},y,x^{(1)}_{\cdot-r}\otimes y)\in M_{d,m}^{\beta}(0,2r)$
\end{enumerate}
\noindent To simplify the proof we will assume $d=m=1$.

\vskip 5pt
\noindent
%%%%%%%%%%%%%%%%%%%%%%%%%%%%%%%%%%%%%%
{\underline{Step 1:}} In order to prove the existence of solution we will use Lemma \ref{puntfix}, a fixed point  argument on $\cc(-r,r,\R_+)$. %We need to define an operator on $\cc(-r,r,\R_+)$.
\vskip 5pt
\noindent
Let us consider the operator 
\[\cl: \cc(-r,r;\R_+)\rightarrow \cc(-r,r;\R_+)\]
such that
\begin{eqnarray*}
\cl(u)(t)&=&\eta(0)+\int_0^t b(s,u)ds+\int_0^t\si(\eta(s-r))dy_s+z(t),\qquad t\in[0,r],\\
\cl(u)(t)&=&\eta(t),\qquad t\in[-r,0].
\end{eqnarray*}
where setting 
\[\xi(t)=\eta(0)+\int_0^t b(s,u)ds+ \int_0^t \si(\eta(s-r))dy_s,\]
then ${\displaystyle z(t)=\max_{s\in[0,t]}(\xi(s))^-}$.
\vskip 5pt
\noindent 
Clearly
$\cl$ is well defined. Let us use the notation $u^*=\cl(u)$.\\
Now, we need to introduce a family of norms in the space $\cc(-r,r;\R_+)$. That is, for any $\la\geq 1$, let us consider
\[\left\|f\right\|_{\infty,\la(-r,r)}:= \sup_{t\in[-r,r]} e^{-\la t}|f(t)|.\]
It is easy to check that these norms are equivalent to $\left\|f\right\|_{\infty(-r,r)}$.
\vskip 5pt
\noindent
Using standard arguments (see for instance \cite{B-R} for similar computations) we obtain that
\begin{eqnarray}
\left\|u^*\right\|_{\infty,\lambda(-r,r)}&\leq&\left\|\eta\right\|_{\infty,\la(-r,0)}+2|\eta(0)|
%\nonumber\\&&
+2 \sup_{t\in[0,r]}e^{-\la t}\left|\int_0^t b(s,u)ds\right|\nonumber\\
&&+2 \sup_{t\in[0,r]}e^{-\la t}\left|\int_0^t \si(\eta(s-r))dy_s\right|.\label{cotaE1}
\end{eqnarray}
We obtain easily (see again  \cite{B-R}) that
{\begin{eqnarray}
\sup_{t\in[0,r]} e^{-\la t}\left|\int_0^t b(s,u)ds\right|&\leq& \frac{L_0}{\la} \left\|u\right\|_{\infty,\la(-r,r)}+\frac{C_\rho}{\la^{1-\rho}}\left\|b_0\right\|_{L^{\rho}}. \label{c1}
\end{eqnarray}}
\vskip 5pt
\noindent 
It only remains the study of the term with the fractional integral. Using the bound appearing on the proof of  Proposition \ref{prop23}, we get for any $\lambda \ge 1$,
\begin{equation}
\begin{array}{l}
\displaystyle\sup_{t\in[0,r]}e^{-\la t}\left|\int_0^t\sigma(\eta({s-r}))dy_s\right|\\
\qquad\quad\displaystyle\leq k|\sigma(\eta({-r}))|\left\|y\right\|_{\beta(0,r)}\sup_{t\in[0,r]}e^{-\la t}t^\beta+\\\qquad\qquad
\displaystyle+k\Phi_{0,r,\beta}(\eta_{\cdot-r},y)\left(\left\|\sigma'\right\|_{\infty}\sup_{t\in[0,r]}e^{-\la t}t^{2\beta}+\left\|\sigma'\right\|_{\gamma}\left\|\eta_{\cdot-r}\right\|^\gamma_{\beta(0,r)}\sup_{t\in[0,r]}e^{-\la t}t^{(\gamma+2)\beta}\right)\\
\qquad\quad\displaystyle\leq k|\sigma(\eta({-r}))|\left\|y\right\|_{\beta(0,r)} \left(\frac{\beta}{\lambda}\right)^\beta e^{-\beta}+\\
\qquad\qquad\displaystyle k\Phi_{0,r,\beta}(\eta_{\cdot-r},y)\left(\left\|\sigma'\right\|_{\infty}\left(\frac{2\beta}{\lambda}\right)^{2\beta} e^{-2\beta}+\left\|\sigma'\right\|_{\gamma}\left\|\eta_{\cdot-r}\right\|^\gamma_{\beta(0,r)}\left(\frac{(\gamma+2)\beta}{\lambda}\right)^{(\gamma+2)\beta} e^{(\gamma+2)\beta}\right)\\
\qquad\quad\displaystyle\leq C_{\beta,\gamma}\frac{1}{\la^{\beta}}\left(|\sigma(\eta({-r}))|\left\|y\right\|_{\beta(0,r)}+\Phi_{0,r,\beta}(\eta_{\cdot-r},y)\left(\left\|\sigma'\right\|_{\infty}+\left\|\sigma'\right\|_{\gamma}\left\|\eta_{\cdot-r}\right\|^\gamma_{\beta(0,r)}\right)\right),
\end{array} \label{a2} 
\end{equation}
where in the last inequality we have used that
\[\sup_{t\in[0,r]} t^\mu e^{-\la t} \leq \left(\frac{\mu}{\la}\right)^\mu e^{-\mu}\]
and $C_{\beta,\gamma}$ is a constant depending on $\beta$ and $\gamma$.
\vskip 5pt
\noindent
So putting together (\ref{cotaE1}), (\ref{c1}) and (\ref{a2}) we have
\[\left\|u^*\right\|_{\infty,\lambda(-r,r)}\leq M_1(\lambda)+M_2(\lambda)\left\|u\right\|_{\infty,\lambda(-r,r)},\]
where
\begin{eqnarray*}
M_1(\lambda)&=&\left\|\eta\right\|_{\infty,\la(-r,0)}+2|\eta(0)|+\frac{2C_\rho}{\la^{1-\rho}}\left\|b_0\right\|_{L^{\rho}}\\
&&+C_{\beta,\gamma}\frac{2}{\la^{\beta}}\left(|\sigma(\eta(-r))|\left\|y\right\|_{\beta(0,r)}+\Phi_{0,r,\beta}(\eta_{\cdot-r},y)\left(\left\|\sigma'\right\|_{\infty}+\left\|\sigma'\right\|_{\gamma}\left\|\eta_{\cdot-r}\right\|^\gamma_{\beta(0,r)}\right)\right),\\
M_2(\la)&=&2L_0\frac{1}{\la}.
\end{eqnarray*}
Now, we can choose $\la=\la_0$ large enough such that $M_2(\la_0) \leq \frac{1}{2}$, Then,  ${\displaystyle\left\|u\right\|_{\infty,\la_0(-r,r)}\leq 2M_1(\la_0)}$ yields that $${\displaystyle\left\|u^*\right\|_{\infty,\la_0(-r,r)}\leq 2M_1(\la_0)}$$ and  ${\displaystyle \cl(B_0)\subseteq B_0}$ for
\[B_0=\left\{u\in \cc(-r,r;\R^d_+); \left\|u\right\|_{\infty,\la_0(-r,r)}\leq 2M_1(\la_0)\right\}.\]
The first hypothesis in Lemma \ref{puntfix} is now satisfied with the metric $\rho_0$ associated to the norm $\left\|\cdot\right\|_{\infty,\la_0(-r,r)}$.
To finish the proof it suffices  to find a metric $\rho_1$ satisfying the second hypothesis in Lemma \ref{puntfix}.
\vskip 5pt \noindent
 Notice first that if $u\in B_0$ then $\left\|u\right\|_{\infty(-r,r)}\leq 2e^{\la_0r}M_1(\la_0):=N_0$.
Consider $u,u'\in B_0$ and $\la\geq 1$. Then
\begin{equation}
 \left\|\cl(u)-\cl(u')\right\|_{\infty,\la(-r,r)}\leq \sup_{t\in[0,r]}e^{-\la t}\left|\xi(t)-\xi'(t)\right|
 +\sup_{t\in[0,r]}e^{-\la t}\left|z(t)-z'(t)\right|.\label{laa}
\end{equation} 
From Lemma \ref{le2} notice that given $t \in [0,r]$ there exists $t_2 \le t$ such that
$$ 
\left|z(t)-z'(t)\right| \le K_l \left|\xi(t_2)-\xi'(t_2)\right| .$$
So
$$ 
e^{-\lambda t}\left|z(t)-z'(t)\right|  \le K_l e^{-\lambda t_2} \left|\xi(t_2)-\xi'(t_2)\right|$$
and it follows easily that
\begin{equation}
\sup_{t\in[0,r]}e^{-\la t}\left|z(t)-z'(t)\right| \le
K_l \sup_{t\in[0,r]}e^{-\la t}\left|\xi(t)-\xi'(t)\right|.\label{lab}\end{equation}
From (\ref{laa}) and (\ref{lab}) we can write
\begin{eqnarray*}
\left\|\cl(u)-\cl(u')\right\|_{\infty,\la(-r,r)}&\leq &(1+K_l)  \sup_{t\in[0,r]}e^{-\la t}\int_0^t\left|b(s,u)-b(s,u')\right|ds\\
&\leq & L_{N_0} (1+K_l) \sup_{t\in[0,r]}e^{-\la t}\int_0^t\sup_{0\leq v\leq s}\left|u(v)-u'(v)\right|ds\\
&\leq & L_{N_0}(1+K_l)\sup_{t\in[0,r]}\int_0^t e^{-\la( t-s)} e^{-\la s} \sup_{-r\leq v\leq s}\left|u(v)-u'(v)\right|ds\\
&\leq & L_{N_0}(1+K_l )\frac{1}{\la}\left\|u-u'\right\|_{\infty,\la(-r,r)} .
\end{eqnarray*}
\vskip 6pt \noindent
So,  choosing $\la=\la_1$ such that ${\displaystyle \frac{L_{N_0}(1+K_l)}{\la_1}\leq \frac{1}{2}}$, the second hypothesis is satisfied for the metric $\rho_1$ associated with the norm $\left\|\cdot\right\|_{\infty,\la_1(-r,r)}$ and $a=\frac12$.
\vskip 7pt
\noindent
%%%%%%%%%%%%%%%%%%%%%%%%%%%%%%%%%%%%%%%
{\underline{Step 2:}} We deal now with the uniqueness problem. \\
Let $x$ and $x'$ be two solutions of (\ref{s1}) in the space $\cc(-r,r;\R_+)$ and choose $N$ large enough such that $\left\|x\right\|_{\infty(-r,r)}\leq N$ and $\left\|x'\right\|_{\infty(-r,r)}\leq N$.
\vskip 5pt
\noindent
For any $t\in[0,r]$,
$$
\sup_{s\in[0,t]} \left|x(s)-x'(s)\right| \le \sup_{s\in[0,t]}\left|\xi(s)-\xi'(s)\right|+ \sup_{s\in[0,t]} \left|z(s)-z'(s)\right|.$$
Moreover, using Lemma \ref{le2} we have
$$
\sup_{s\in[0,t]} \left|z(s)-z'(s)\right|\leq K_l \sup_{s\in[0,t]} \left|\xi(t)-\xi'(t)\right|.$$
So, putting together the last two inequalities we get that
\begin{eqnarray*}
\sup_{s\in[0,t]} \left|x(s)-x'(s)\right|&\le & (1+K_l)  \sup_{s\in[0,t]} \left|\xi(s)-\xi'(s)\right| \\&\leq&
(1+K_l)  \sup_{s\in[0,t]}
\left|\int_0^s \left(b(\tau,x)-b(\tau,x')\right)d\tau\right|\\&\leq&
(1+K_l) L_N\sup_{s\in[0,t]}\left|\int_0^s\sup_{0\leq v\leq \tau}|x(v)-x'(v)|d\tau\right|\\&\leq& L_N(1+K_l)\int_0^t\sup_{v\in[0, \tau]}|x(v)-x'(v)|d\tau.
\end{eqnarray*}
Applying now Gronwall's inequality, we have that for all $t\in[0,r]$
\[\sup_{s\in[0,t]}|x(s)-x'(s)|= 0.\]
So
\[\left\|x-x'\right\|_{\infty(-r,r)}= 0 \]
and the uniqueness has been proved.
\vskip 5pt
\noindent
%%%%%%%%%%%%%%%%%%%%%%%%%
{\underline{Step 3:}} 
We have to prove that  $(x^{(1)}_{\cdot-r},y,x^{(1)}_{\cdot-r}\otimes y)\in M_{1,1}^{\beta}(0,2r)$.\\

We have to check the three conditions appearing in Definition \ref{d2.1}:
\begin{enumerate}
\item $y:[0,2r]\rightarrow\R$ is $\beta$-H\"older continuous. This condition is one of the hypothesis of our theorem.
\item $x^{(1)}_{\cdot-r}:[0,2r]\rightarrow \R$ is $\beta$-H\"older continuous.\\
We can write that
\begin{eqnarray*}
\left\|x^{(1)}_{\cdot-r}\right\|_{\beta(0,2r)}&=&\left\|x^{(1)}\right\|_{\beta(-r,r)}=\sup_{-r\leq v\leq w\leq r} \frac{|x^{(1)}(w)-x^{(1)}(v)|}{(w-v)^\beta}\\
&\leq& \sup_{-r\leq v\leq w<0} \frac{|\eta(w)-\eta(v)|}{(w-v)^\beta}+\sup_{\substack{-r\leq v\leq 0\\ 0\leq w\leq r}}\frac{|x^{(1)}(w)-\eta(v)|}{(w-v)^\beta}\\
&&+ \sup_{0\leq v\leq w\leq r} \frac{|x^{(1)}(w)-x^{(1)}(v)|}{(w-v)^\beta}.
\end{eqnarray*}
Note that 
\[\frac{|x^{(1)}(w)-\eta(v)|}{(w-v)^\beta}\leq\frac{|x^{(1)}(w)-\eta(0)|}{(w-0)^\beta}+\frac{|\eta(0)-\eta(v)|}{(0-v)^\beta}. \]
So
\begin{equation}
\left\|x^{(1)}_{\cdot-r}\right\|_{\beta(0,2r)}\leq 2\left\|\eta\right\|_{\beta(-r,0)}+2\left\|x^{(1)}\right\|_{\beta(0,r)}. \label{a3}
\end{equation}

%Comencem estudiant el terme $A_1$. Per a $w\in[0,r]$ i $i\in \left\{1,\ldots,d\right\}$ denotem per
%\[w_0^i=\inf\left\{u;z^i(u)=z^i(w)\right\}\vee 0,\]
%aleshores es compleix que
%\[|z^i(w)|=|\xi^i(w_0^i)|\leq \sup_{0\leq s\leq w^i_0}|\xi^i(s)|.\]
%I per tant podem acotar el terme $A_1$ de la seg�ent manera,
%\begin{eqnarray*}
%A_1&=&\sup_{\substack{-r\leq v\leq 0\\ 0\leq w\leq r}}\left(\left|\eta(0)+\int_0^w b(s,x)ds+\int_0^w \sigma(\eta_{s-r})dy_s+z(w)-\eta(v)\right|\right)(w-v)^{-\beta}\\
%&\leq&(d+1)\sup_{\substack{-r\leq v\leq 0\\ 0\leq w\leq r}}\frac{|\eta(0)-\eta(v)|}{(w-v)^\beta}+(d+1)\sup_{\substack{-r\leq v\leq 0\\ 0\leq w\leq r}}\frac{\left|\int_0^wb(s,x)ds\right|}{(w-v)^\beta}\\&&+(d+1)\sup_{\substack{-r\leq v\leq 0\\ 0\leq w\leq r}}\frac{\int_0^w \sigma(\eta_{s-r})dy_s}{(w-v)^\beta}\\
%&\leq&(d+1)\left\|\eta\right\|_{\beta(-r,0)}+(d+1)\left\|\int b(s,x)ds\right\|_{\beta(0,r)}+(d+1)\left\|\int \sigma(\eta_{s-r})dy_s\right\|_{\beta(0,r)},
%\end{eqnarray*}
%on hem utilitzant que $(w-v)^\beta\leq(w-0)^\beta$.
%Pel que fa al segon terme tenim que 

Moreover
\begin{equation}
\left\|x^{(1)}\right\|_{\beta(0,r)} \leq \left\|\int_0^\cdot b(s,x^{(1)})ds\right\|_{\beta(0,r)}+\left\|\int_0^\cdot \sigma(\eta(s-r))dy_s\right\|_{\beta(0,r)}+\left\|z^{(1)}\right\|_{\beta(0,r)}. \label{a4}
\end{equation}
Using Lemma \ref{le2} we also get that
\begin{equation}
\left\|z^{(1)}\right\|_{\beta(0,r)}\leq \left\|\xi^{(1)}\right\|_{\beta(0,r)}. \label{a5}
\end{equation}

%Pel que fa la norma $\left\|z(\cdot)\right\|_{\beta(0,r)}$, sabem que si $z$ �s creixent en $t$ aleshores $z^i(t)=-\xi^i(t)$, i tamb� que, per qualsevol $s$ sempre �s cert que $z^i(s)\geq -\xi^i(s)$, per tant 
%\begin{equation*}
%|z^i(t)-z^i(s)|=z^i(t)-z^i(s)\leq -\xi^i(t)+ \xi^i(s)\leq |\xi^i(t)-\xi^i(s)|.
%\end{equation*}
%amb aquest resultat i si per a $w\in[0,r]$ i $i\in \left\{1,\ldots,d\right\}$ denotem per
%\[w_0^i=\inf\left\{u;z^i(u)=z^i(w)\right\}\vee 0,\]
%�leshores tamb� tenim que
%\[\frac{|z^i(w)-z^i(v)|}{(w-v)^\beta}=\frac{|z^i(w_0^i)-z^i(v)|}{(w-v)^\beta}\leq \frac{|z^i(w_0^i)-z^i(v)|}{(w_0^i-v)^\beta}\leq \frac{|\xi^i(w_0^i)-\xi^i(v)|}{(w_0^i-v)^\beta}.\]
%Per tant,
%\[\left\|z(\cdot)\right\|_{\beta(0,r)}\leq d\left(\left\|\int_0^\cdot b(s,x)ds\right\|_{\beta(0,r)}+\left\|\int_0^\cdot \sigma(\eta_{s-r})dy_s\right\|_{\beta(0,r)}\right).\]
%I llavors 
%\[A\leq (d+1)\left(\left\|\int_0^\cdot b(s,x)ds\right\|_{\beta(0,r)}+\left\|\int_0^\cdot \sigma(\eta_{s-r})dy_s\right\|_{\beta(0,r)}\right).\]
%En els treballs [F-R] i [H-N] trobem les seg�ents cotes,
%\begin{eqnarray*}
%\left\|\int_0^\cdot b(s,x)ds\right\|_{\beta(0,T)}&\leq& d^{(1)}(1+\left\|x\right\|_{\infty(r)}),\\
%\left\|\int_0^\cdot \si(x_s)ds\right\|_{\beta(a,b)}&\leq& k|\si(x_a)|\left\|y\right\|_{\beta(a,b)}+k\Phi_{a,b,\beta}(x,y)\\
%&&\times\left(\left\|\si'\right\|_\infty+\left\|\si'\right\|_\la\left\|x\right\|^\la_{\beta(a,b)}(b-a)^{\la\beta}\right)(b-a)^\beta.
%\end{eqnarray*}
Furthermore,  putting together (\ref{a3}), (\ref{a4}) and (\ref{a5}) and using again Proposition  \ref{prop23}  we obtain that
\begin{eqnarray*}
\left\|x^{(1)}\right\|_{\beta(-r,r)}&\leq&4\left\|\eta\right\|_{\beta(-r,0)}+4\left\|\int_0^\cdot b(s,x^{(1)} )ds\right\|_{\beta(0,r)}
%\\&&
+4\left\|\int_0^\cdot \sigma(\eta (s-r))dy_s\right\|_{\beta(0,r)}\\
&\leq&4\left\|\eta\right\|_{\beta(-r,0)}+4 C\left(1+\left\|x^{(1)}\right\|_{\infty(-r,r)}\right)
%\\&&
+4(k|\sigma(\eta (-r))|\left\|y\right\|_{\beta(0,r)}\\
&&+\Phi_{0,r,\beta}(\eta_{\cdot-r},y)(\left\|\sigma'\right\|_\infty+\left\|\sigma'\right\|_\gamma\left\|\eta_{\cdot-r}\right\|^\gamma_{\beta(0,r)}r^{\gamma\beta}))r^\beta.
\end{eqnarray*}
So we can conclude that $x^{(1)}_{\cdot-r}$ is $\beta$-H\"older continuous.

\item Let us define $(x^{(1)}_{\cdot-r}\otimes y)_{s,t}$ for $s,t\in \Delta_{2r}$. For completeness, we will give this definition for any dimensions $d$ and $m$, unless we will still  consider  $d=m=1$ in the proofs. For any $k \in \{1,\cdots,d\}$ and $l \in \{1,\cdots,m\}$, set:
\begin{itemize}
\item
if $s,\,t\in[0,r]$, 
\[(x^{(1)}_{\cdot-r}\otimes y)_{s,t}^{k,l}=(\eta_{\cdot-r}\otimes y)_{s,t}^{k,l},\]
\item if $s,\,t\in [r,2r]$, set
\begin{eqnarray*}
(x^{(1)}_{\cdot-r}\otimes y)_{s,t}^{k,l}&=&\int_s^t (y^l(t)-y^l(v))b^k(v-r,x^{(1)})dv+\sum_{j=1}^m \int_{s}^{t} \si^k_j(\eta(v-2r))d(y_{\cdot-r}\otimes y)_{\cdot, t}^{j,l}(v)\\
&&+\int_s^t (y^l(t)-y^l(v))d(z^{(1)})_{v-r}^k,
\end{eqnarray*}
\item if $s\in[0,r]$ and $t\in[r,2r]$, 
\[(x^{(1)}_{\cdot-r}\otimes y)_{s,t}^{k,l}=(\eta_{\cdot-r}\otimes y)_{s,r}^{k,l}+(x^{(1)}_{\cdot-r}\otimes y)_{r,t}^{k,l}+(\eta^k(0)-\eta^k(s-r))\otimes (y^l(t)-y^l(r)).\]
\end{itemize}
\vskip 5pt
\noindent
Let us check that the multiplicative property (let us recall that we consider again $d=m=1$ for simplicity) is satisfied, that is, 
for any $ 0 \le s \le u \le t \le 2r$ it holds that
\begin{equation}
(x^{(1)}_{\cdot-r}\otimes y)_{s,u}+(x^{(1)}_{\cdot-r}\otimes y)_{u,t}+(x^{(1)}(u-r)-x^{(1)}(s-r))\otimes (y(t)-y(u))=(x^{(1)}_{\cdot-r}\otimes y)_{s,t}. \label{s2}
\end{equation}

We have to distinguish several cases:
\begin{itemize}
%%%%%%%%%%%%%%%%%%%%%%%%
\item[a)] Case $0\leq s\leq u\leq t\leq r$.
 
Since on $\Delta_r$ it holds that
$$(x^{(1)}_{\cdot-r},y,x^{(1)}_{\cdot-r}\otimes y)=(\eta_{\cdot-r},y,\eta_{\cdot-r}\otimes y),$$ 
the multiplicative property follows from the fact that  we are assuming  that $(\eta_{\cdot-r},y,\eta_{\cdot-r}\otimes y)$  
 is a  $\beta-$ H\"older continuous functional.

%%%%%%%%%%%%%%%%%%%%%%%%%%%%%%%%%%%%%%%%%%%%%
\item[b)] Case  $r\leq s\leq u\leq t\leq 2r$. 

Notice first that,
\begin{equation*}
\begin{array}{l}
\displaystyle (x^{(1)}_{\cdot-r}\otimes y)_{s,u}+(x^{(1)}_{\cdot-r}\otimes y)_{u,t}=\int_s^u (y(u)-y(v))b(v-r,x^{(1)})dv\\[4mm]
\displaystyle \qquad\quad+\int_{s}^{u} \sigma(\eta(v-2r))d(y_{\cdot-r}\otimes y)_{\cdot,u}(v)+\int_s^u(y(u)-y(v))dz^{(1)}_{v-r}\\[4mm]
\displaystyle \qquad\quad+\int_u^t (y(t)-y(v))b({v-r},x^{(1)})dv+\int_{u}^{t} \sigma(\eta(v-2r))d(y_{\cdot-r}\otimes y)_{\cdot,t}(v)\\[4mm]
\displaystyle \qquad\quad+\int_u^t(y(t)-y(v))dz^{(1)}_{v-r}\\[4mm]
\displaystyle \qquad =\int_s^t (y(t)-y(v))b({v-r},x^{(1)})dv+\int_{s}^{t} \sigma(\eta(v-2r))d(y_{\cdot-r}\otimes y)_{\cdot,t}(v)\\[4mm]
\displaystyle \qquad\quad +\int_s^t(y(t)-y(v))dz^{(1)}_{v-r}\\[4mm]
\displaystyle \qquad\quad+(y(u)-y(t))\left(\int_s^u b({v-r},x^{(1)})dv+z^{(1)}(u-r)-z^{(1)}(s-r)\right)\\[4mm]
\displaystyle \qquad\quad+\int_{s}^{u} \si(\eta(v-2r))(d(y_{\cdot-r}\otimes y)_{\cdot,u}(v)-d(y_{\cdot-r}\otimes y)_{\cdot,t}(v)).
\end{array}
\end{equation*}
So
\begin{equation}
\begin{array}{l}
\displaystyle (x^{(1)}_{\cdot-r}\otimes y)_{s,u}+(x^{(1)}_{\cdot-r}\otimes y)_{u,t}\\[4mm]
\displaystyle \qquad =(x^{(1)}_{\cdot-r}\otimes y)_{s,t}+(y(u)-y(t))\left(\int_s^u b({v-r},x^{(1)})dv+z^{(1)}(u-r)-z^{(1)}(s-r)\right)
\\[4mm]
\displaystyle \qquad\quad+\int_{s}^{u} \si(\eta(v-2r))(d(y_{\cdot-r}\otimes y)_{\cdot,u}(v)-d(y_{\cdot-r}\otimes y)_{\cdot,t}(v))
\end{array}\label{xxa}
\end{equation}
On the other hand, from Definition \ref{d2.1} we obtain that
\begin{equation}
\int_{s}^{u} \si(\eta(v-2r))(d(y_{\cdot-r}\otimes y)_{\cdot,u}(v)-d(y_{\cdot-r}\otimes y)_{\cdot,t}(v))=(y(u)-y(t))\int_{s}^{u} \si(\eta(v-2r))dy_{v-r}.
\label{xxb}
\end{equation}
Finally, using that
\begin{eqnarray*}
\int_s^u b({v-r},x^{(1)})dv&=&\int_{s-r}^{u-r} b({v},x^{(1)})dv,\\
\int_{s}^{u} \si(\eta(v-2r))dy_{v-r}&=&\int_{s-r}^{u-r} \si(\eta(v-r))dy_{v},
\end{eqnarray*}
and putting together (\ref{xxa}) and (\ref{xxb}) we get the multiplicative property (\ref{s2}).

%%%%%%%%%%%%%%%%%%%%%%%%%%%%%%%%%%%%%%%%%%%%%%%%%%%%%%%%
\item[c)] Case  $0\leq s\leq r$ and $r\leq u\leq t\leq 2r$.\\
Notice first that from the definition of $(x^{(1)}_{\cdot-r}\otimes y)$  it follows that
\begin{equation}
(x^{(1)}_{\cdot-r}\otimes y)_{s,u}=(\eta_{\cdot-r}\otimes y)_{s,r}+(x^{(1)}_{\cdot-r}\otimes y)_{r,u}+(\eta(0)-\eta(s-r))\otimes (y(u)-y(r)).
\label{bba}
\end{equation}
On the other hand, we have seen in the case b) (choosing $s=r$) that
\begin{equation}
(x^{(1)}_{\cdot-r}\otimes y)_{r,t}=(x^{(1)}_{\cdot-r}\otimes y)_{r,u}+(x^{(1)}_{\cdot-r}\otimes y)_{u,t}+(x^{(1)}(u-r)-\eta(0))\otimes (y(t)-y(u)). 
\label{bbb}
\end{equation}
So, putting together (\ref{bba}) and (\ref{bbb}) we can write
\begin{equation*}
\begin{array}{l}
\displaystyle (x^{(1)}_{\cdot-r}\otimes y)_{s,u}+(x^{(1)}_{\cdot-r}\otimes y)_{u,t}+(x^{(1)}(u-r)-\eta(s-r))\otimes (y(t)-y(u))\\[4mm]
\displaystyle \qquad =(\eta_{\cdot-r}\otimes y)_{s,r}+(x^{(1)}_{\cdot-r}\otimes y)_{r,t}+(\eta_0-\eta_{s-r})\otimes (y_t-y_r)\\[4mm]
%\displaystyle \qquad\quad +(x^1_{\cdot-r}\otimes y)_{u,t}+(x^1(u-r)-\eta(s-r))\otimes (y(t)-y(u))\\[4mm]
%\displaystyle \qquad =(\eta_{\cdot-r}\otimes y)_{s,r}+(x^1_{\cdot-r}\otimes y)_{r,u}+(x^1_{\cdot-r}\otimes y)_{u,t}+(\eta_0-\eta_{s-r})\otimes (y_
%t-y_r)\\[4mm]
\displaystyle \qquad=(x^{(1)}_{\cdot-r}\otimes y)_{s,t},
\end{array}
\end{equation*}
where the last equality follows for the definition of $(x^{(1)}_{\cdot-r}\otimes y)$. The proof of this case is now finished.
%%%%%%%%%
\item[d)] Case $0\leq s\leq u\leq r$ and $r\leq t\leq 2r$.

This case can be done following the same ideas that the case c).
%On the one hand,
%\[(x^1_{\cdot-r}\otimes y)_{s,t}=(\eta_{\cdot-r}\otimes y)_{s,r}+(x^1_{\cdot-r}\otimes y)_{r,t}+(\eta_{0}-\eta_{s-r})\otimes (y_t-y_r).\]
%On the other hand,
%\begin{equation*}
%\begin{array}{l}
%\displaystyle (\eta_{\cdot-r}\otimes y)_{s,u}+(x^1_{\cdot-r}\otimes y)_{u,t}+(\eta_{u-r}-\eta_{s-r})\otimes (y_t-y_u)\\[4mm]
%\displaystyle \qquad = (\eta_{\cdot-r}\otimes y)_{s,u}+(\eta_{\cdot-r}\otimes y)_{u,r}+(x^1_{\cdot-r}\otimes y)_{r,t}+(\eta_{0}-\eta_{u-r})\otimes %(y_t-y_r)\\[4mm]
%\displaystyle \qquad\quad+(\eta_{u-r}-\eta_{s-r})\otimes (y_t-y_u)\\[4mm]
%\displaystyle \qquad =(\eta_{\cdot-r}\otimes y)_{s,u}+(\eta_{\cdot-r}\otimes y)_{u,r}+(\eta_{u-r}-\eta_{s-r})\otimes (y_r-y_u)+(x^1_{\cdot-r}\otimes %y)_{r,t}\\[4mm]
%\displaystyle \qquad\quad+(\eta_{0}-\eta_{s-r})\otimes (y_t-y_r)\\[4mm]
%\displaystyle \qquad =(\eta_{\cdot-r}\otimes y)_{s,r}+(x^1_{\cdot-r}\otimes y)_{r,t}+(\eta_{0}-\eta_{s-r})\otimes (y_t-y_r).
%\end{array}
%\end{equation*}
\end{itemize}
%%%%%%%%%%%%%%%%%%%%%%%%%%%%%%%%%%%%%%%%

\item Now only remains to prove that $|(x^{(1)}_{\cdot-r}\otimes y)_{s,t}|\leq C|t-s|^{2\beta}$. We will distinguish again three cases:\\

\begin{enumerate}

\item Assume that $s,t\in[r,2r]$. Then
\begin{eqnarray*}
|(x^{(1)}_{\cdot-r}\otimes y)_{s,t}|&\leq& \left|\int_s^t(y(t)-y(v))b({v-r},x^{(1)})dv\right|+\left|\int_{s}^{t} \sigma(\eta(v-2r))d(y_{\cdot-r}\otimes y)_{\cdot,t}(v)\right|\\ &&+\left|\int_s^t(y(t)-y(v))dz^{(1)}_{v-r}\right|.
\end{eqnarray*}
Since $y$ is $\beta-$H\"older continuous function, we have that
\[\left|\int_s^t (y(t)-y(v))dz^{(1)}_{v-r}\right|\leq K|t-s|^\beta|z^{(1)}(t-r)-z^{(1)}(s-r)|,\]
for a constant $K$. 
Then, using Lemma \ref{le2} we get
\begin{equation}
\left|\int_s^t (y(t)-y(v))dz^{(1)}_{v-r}\right|\leq K|t-s|^{2\beta}. \label{s3}
\end{equation}
On the other hand, using the hypothesis on $b$ we have
\begin{eqnarray}
& &\left|\int_s^t(y(t)-y(v))b(v-r,x^{(1)})dv\right| \leq K|t-s|^\beta\left|\int_s^t(L_0\sup_{-r\leq u\leq v-r}|x^{(1)}(u)|+b_0(v))dv\right|\nonumber\\
& & \qquad  \leq K|t-s|^{\beta+1}\left\|x^{(1)}\right\|_{\infty(-r,r)}+|t-s|^{\beta+1-\frac{1}{\rho}}\left\|b_0\right\|_{L^\rho}. \label{s4}
\end{eqnarray}
Finnally using Proposition  \ref{prop24} we get
\begin{equation}
\begin{array}{l}
\displaystyle 
\left|\int_{s}^{t} \sigma(\eta(v-2r))d(y_{\cdot-r}\otimes y)_{\cdot,t}(v)\right|\leq k|\sigma(\eta(s-2r))|\Phi_{s,t,\beta}(y_{\cdot-r},y)(t-s)^{2\beta}\\[4mm]
\displaystyle \qquad\quad+k\left(\left\|\sigma'\right\|_\infty+\left\|\sigma'\right\|_\gamma\left\|\eta_{.-2r}\right\|^\gamma_{\beta(s,t)}(t-s)^{\gamma\beta}\right)\Phi_{s,t,\beta}(\eta_{\cdot-2r},y_{\cdot-r},y)(t-s)^{3\beta},
\end{array} \label{s5}
\end{equation}
where 
\begin{eqnarray*}
\Phi_{a,b,\beta}(x,y,z)&=&\left\|y\right\|_{\beta(a,b)}\left\|z\right\|_{\beta(a,b)} \left\|x\right\|_{\beta(a,b)} \\
&& +\left\|z\right\|_{\beta(a,b)}\left\|x\otimes y\right\|_{
2\beta(a,b)}+\left\|x\right\|_{\beta(a,b)}\left\|y\otimes z\right\|_{2\beta(a,b)}.
\end{eqnarray*}
\vskip 5pt
\noindent
Now putting together (\ref{s3}), (\ref{s4}) and (\ref{s5}) we finish the proof.
\vskip 5pt
\noindent
\item If $s\in[0,r]$ and $t\in[r,2r]$,
\begin{eqnarray*}
|(x^{(1)}_{\cdot-r}\otimes y)_{s,t}|&\leq&|(\eta_{\cdot-r}\otimes y)_{s,r}|+|(x^{(1)}_{\cdot-r}\otimes y)_{r,t}|+|(\eta(0)-\eta(s-r))\otimes (y(t)-y(r))|\\
&\leq& K|r-s|^{2\beta}+K|t-r|^{2\beta}+K|s-r|^\beta|t-r|^\beta\leq K|t-s|^{2\beta}.
\end{eqnarray*}
\item  If $s,\,t\in[0,r]$ then $x^{(1)}_{\cdot-r}=\eta_{\cdot-r}$ and the result is already true.
\end{enumerate}
\end{enumerate}
\hfill $\Box$

%%%%%%%%%%%%%%%%%%%%%%%%%%%%%%%%%%%%%%%%%%%%%%%%%%%%%%%%%%%%%%%%%%%%%%%%%%%
\renewcommand{\theequation}{5.\arabic{equation}} \setcounter{equation}{0}

\section{Boundedness for deterministic integral equations}

The aim of this section is the proof of Theorem \ref{tt2}. For simplicity let us assume $T=Mr$.

\vskip 5pt
\noindent
{\bf Proof of Theorem \ref{tt2}:}
The proof will be done in several steps.
\vskip 5pt
\noindent
%%%%%%%%%%%%%%%%%%%%%%%%%
{\underline{Step 1:}}  Assuming that $(x,y,x\otimes y)\in M^\beta_{d,m}(0,T)$, let us  define 
$\left(x_{\cdot-r}\otimes y_{\cdot-r}\right)_{s,t}$. Set
\begin{equation}\label{e5.1biss}\left(x_{\cdot-r}\otimes y_{\cdot-r}\right)_{s,t}:=\left(x\otimes y\right)_{s-r,t-r}.\end{equation}
It clearly t belongs to $M^\beta_{d,m}(r,T)$. Notice that the functions $x_{\cdot-r}$ and $y_{\cdot-r}$ are $\beta$-H\"older
continuous and $x_{\cdot-r}\otimes y_{\cdot-r}$ is a continuous functions satisfying the multiplicative property. Indeed, we have, for 
$s\leq u\leq t$, 
\begin{equation*}\begin{array}{l}
\displaystyle (x_{\cdot-r}\otimes y_{\cdot-r})_{s,u}+(x_{\cdot-r}\otimes y_{\cdot-r})_{u,t}+(x(u-r)-x(s-r))\otimes
(y(t-r)-y(u-r))\\
\qquad \displaystyle 
=(x\otimes y)_{s-r,u-r}+(x\otimes y)_{u-r,t-r}+(x(u-r)-x(s-r))\otimes
(y(t-r)-y(u-r))\\
\qquad \displaystyle =(x\otimes y)_{s-r,t-r} =(x_{\cdot-r}\otimes y_{\cdot-r})_{s,t}.
\end{array}\end{equation*}
Finally, we also have that,
for all $(s,t)\in \{(s,t): r\le s<t\le T\}$, 
\begin{equation*}
\left| (x_{\cdot-r}\otimes y_{\cdot-r})_{s,t}\right|= \left| (x\otimes y)_{s-r,t-r}\right|
\leq c|t-s|^{2\beta }.
\end{equation*}
\vskip 5pt
\noindent
%%%%%%%%%%%%%%%%%%%%%%%%%
{\underline{Step 2:}}  Set, for any $s, t \in [nr,(n+1)r]$, $n=0,\dots, (M-1)$,
\begin{eqnarray*}
J_1^n(x,y,x\otimes y) (t)&=&\eta(0)+\int_0^t b(s,x)ds+\int_0^t\sigma(x(s-r))dy_s+ z(t),\\
J_2^n(x,y,x\otimes y )(s,t)&=&
\left\{\begin{array}{ll}
(\eta_{\cdot-r}\otimes y)_{s,t}, & {\rm if}\ s,t\in [0,r] \, {\rm and} \, n=1,\\
(x_{\cdot-r}\otimes y)_{s,t}, & {\rm otherwise}, 
\end{array}\right.
\end{eqnarray*}
where
\begin{eqnarray*}
(x_{\cdot-r}\otimes y)_{s,t}&=&\int_s^t (y(t)-y(u)) b(u-r,x)du+\int_{s}^{t}\sigma(x(u-2r))d(y_{\cdot-r}\otimes y)_{\cdot,t}(u)\\
&&
+ \int_s^t(y(t)-y(u))dz_{u-r},\end{eqnarray*}
with $z$ given in (\ref{det}). 

\noindent Set
 \begin{displaymath}
 \Delta_y= \left(\left\|\eta\right\|_{\beta}
+\left\|\eta_{\cdot-r}\otimes y\right\|_{2\beta}
+(d^\frac12 +1)(1+3k) \mu \left[\left\|y\right\|_{\beta}
+\left\|y\right\|_{\beta}^2+ \left\|y_{\cdot-r}\otimes y\right\|_{2\beta}\right]\right)^{-\frac{1}{\beta}},
  \end{displaymath}
  with $k$ depending only on $\beta$ and $\gamma$.
If we consider $s,t$ such that
$$0<t-s\le \Delta_y \wedge 1,$$
then we will have that
\begin{eqnarray}
(t-s)^\beta&\le &\Delta_y^\beta\le \frac{1}{\left\|\eta\right\|_{\beta}},\label{e5.7}\\[2mm]
(t-s)^\beta&\le &\Delta_y^\beta\le \frac{1}{\left\|\eta_{\cdot-r}\otimes y\right\|_{2\beta}},\label{e5.8}\\[2mm]
(t-s)^\beta&\le &\Delta_y^\beta\le \frac{1}{(d^\frac12 +1)(1+3k)\ \mu\left(1 +\left\|y\right\|_{\beta}
+\left\|y\right\|_{\beta}^2+ \left\|y_{\cdot-r}\otimes y\right\|_{2\beta}\right)}.\label{e5.9}
\end{eqnarray}
\vskip 5pt
%%%%%%%%%%%%%%%%%%%%%%%%%
\noindent We will use an induction argument  to prove that for any $n$ and for all 
$s,t$  such that $nr \le s < t \le (n+1)r$ with $0<t-s\le \Delta_y \wedge 1,$ it hols that
$$\Vert J_1^n \Vert_{\beta(s,t)} = \Vert x \Vert_{\beta(s,t)} \le \mu (d^\frac12 +1)\left[1 +2k+ 3k \left\|y\right\|_{\beta}\right].$$

\vskip 5pt
%%%%%%%%%%%%%%%%%%%%%%%%%
\noindent
Assuming this last inequality, the proof of (\ref{e5.1}) is standard. Indeed, notice that
it follows easily, that  for any  $0 \le s<t \le T$ with  $t-s\le \Delta_y$, 
$$ \Vert x \Vert_{\beta(s,t)} \le 2 \mu (d^\frac12 +1) \left[1+2k+3k \left\|y\right\|_{\beta}\right].$$
Then,  if $t-s\le \Delta_y$, $0\le s<t$, we have
\begin{equation}\label{e5.16}
\sup_{r\in[s,t]} \left|x(r)\right|\le \left|x(s)\right| + (t-s)^\beta \left\|x\right\|_{\beta(s,t)}\le  \left|x(s)\right| +2,
\end{equation}
and, particularly,  
$$\sup_{r\in[0,\Delta_y]} \left|x(r)\right| \le  \eta(0) +2.$$
Now we divide the interval $[0,T]$ into $n=\left[T\Delta_y^{-1}\right]+1$ intervals of length $\Delta_y$, where $[a]$ denotes the largest integer bounded by $a$. Then, applying (\ref{e5.16}) on the intervals $\left[\Delta_y,2\Delta_y\right], \dots, \left[(n-1)\Delta_y,n\Delta_y\right]$ and the previous inequality, we obtain
\begin{displaymath}
\sup_{r\in[0,T]} \left|x(r)\right|\le \eta(0) +2n \le 2+\eta(0) + 2\left[T\Delta_y^{-1}\right],
\end{displaymath}
and we can conclude that the estimate (\ref{e5.1}) is true.

\vskip 5pt
%%%%%%%%%%%%%%%%%%%%%%%%%
\noindent Let us come back to check our induction argument to finish the proof.

\vskip 5pt
%%%%%%%%%%%%%%%%%%%%%%%%%
\noindent {\underline{Step 2.1:}}
Assume  $s,t \in [0,r]$. 
On the one hand, 
\begin{equation*}%\label{e5.10}
\left\|J_2^0\right\|_{2\beta (s,t)}=\left\|\eta_{\cdot-r}\otimes y\right\|_{2\beta (s,t)}.
\end{equation*}
On the other hand, using Lemma \ref{le2} and  Proposition  \ref{prop23}, we obtain that
\begin{eqnarray*}
\left\|J_1^0\right\|_{\beta (s,t)}&\le & (d^\frac12 +1)\bigg[\|b\|_\infty + k\Big(\left\|\sigma\right\|_{\infty} \left\|y\right\|_\beta+ \left[\left\|\eta_{\cdot-r}\otimes y\right\|_{2\beta (s,t)}+\|\eta_{\cdot-r}\|_{\beta(s,t)}\left\|y\right\|_\beta\right]\\
&&\times\left[\left\| \sigma' \right\|_{\infty}+
\left\| \sigma' \right\|_{\gamma} \left\| \eta_{\cdot-r} \right\|_{\beta(s,t)}^\gamma(t-s)^{\gamma \beta}\right](t-s)^\beta\Big)\bigg].
\end{eqnarray*}
By means of  (\ref{e5.7}) and (\ref{e5.8}), we have
\begin{eqnarray*}
\left\| \eta_{\cdot-r} \right\|_{\beta(s,t)}(t-s)^\beta&\le &1,\\
 \left\|\eta_{\cdot-r}\otimes y\right\|_{2\beta (s,t)} (t-s)^\beta&\le &1.\end{eqnarray*}
Then we conclude that
\begin{equation}\label{e5.11}
\left\|J_1^0\right\|_{\beta (s,t)}\le  \mu (d^\frac12 +1)\left[1 +2k+ 3k \left\|y\right\|_{\beta}\right].
\end{equation}
\vskip 5pt
\noindent
%%%%%%%%%%%%%%%%%%%%%%%%%
{\underline{Step 2.2:}} Since the first step does not follow the general case, in order to prove our induction we have also to consider the next interval.

\noindent So, assume now $s,t\in [r,2r]$. 
We first need to deal with $\left\|J_2^1\right\|_{2\beta (s,t)}=\left\|x_{\cdot-r}\otimes y\right\|_{2\beta (s,t)}$.
Applying Proposition  \ref{prop24} and dealing with the integrals with $b$ and $z$, we get
\begin{eqnarray*}
\left\|J_2^1\right\|_{2\beta (s,t)}&\le & \left\|y\right\|_\beta \left\|z\right\|_{\beta(s-r,t-r)}
+\mu \bigg[ \left\|y\right\|_\beta  +k\Big(\left\|y_{\cdot-r}\otimes y\right\|_{2\beta}
+\left\|y\right\|_{\beta}^2\\&&+\left[1+ \left\| x_{.-2r} \right\|_{\beta(s,t)}^\gamma(t-s)^{\gamma \beta}\right]
 \Big[\left\| x_{.-2r} \right\|_{\beta(s,t)}\left\|y\right\|_{\beta}^2\\
 &&+\left\|x_{.-2r}\otimes y_{\cdot-r}\right\|_{2\beta (s,t)}\left\|y\right\|_{\beta}+\left\| x_{.-2r} \right\|_{\beta(s,t)}
\left\|y_{\cdot-r}\otimes y\right\|_{2\beta}\Big](t-s)^\beta\Big)\bigg].
\end{eqnarray*}
The same computations given to bound $J_1^0$ in (\ref{e5.11}) and Lemma \ref{le2} yields that
\begin{displaymath}
\left\|z\right\|_{\beta(s-r,t-r)}\le d^\frac12  \mu \left[1 +2k+3k \left\|y\right\|_{\beta}\right].
\end{displaymath}
As before, the estimate (\ref{e5.7}) gives
\begin{displaymath}
\left\| x_{.-2r} \right\|_{\beta(s,t)}(t-s)^\beta=\left\| \eta_{.-2r} \right\|_{\beta(s,t)}(t-s)^\beta\le 1,
\end{displaymath}
and (\ref{e5.1biss}) and (\ref{e5.8}) imply
$$
\left\| x_{.-2r}\otimes y_{\cdot-r} \right\|_{2\beta(s,t)}(t-s)^\beta = \left\| \eta_{\cdot-r}\otimes y  \right\|_{2\beta(s-r,t-r)}(t-s)^\beta
%\\
%&=&\left\| J_2^0 \right\|_{2\beta(s-r,t-r)}(t-s)^\beta\\
 \le 1.
$$
So, using all these inequalities, we can bound $J_2^1$ as follows:
\begin{eqnarray}\label{e5.12}
\left\|J_2^1\right\|_{2\beta (s,t)}&\le & (d^\frac12 +1)
(1+2k) \mu \left\|y\right\|_{\beta}+(d^\frac12 +1) 3 k \mu \left\|y\right\|_{\beta}^2+
3 k\mu \left\|y_{\cdot-r}\otimes y\right\|_{2\beta}\nonumber\\
&\le & 
 (d^\frac12 +1)
(1+3k) \mu \left(\left\|y\right\|_{\beta}+
\left\|y\right\|_{\beta}^2+ \left\|y_{\cdot-r}\otimes y\right\|_{2\beta}\right).
\end{eqnarray}
Proposition  \ref{prop23}  again yields
\begin{eqnarray*}
\left\|J_1^1\right\|_{\beta (s,t)}&\le &  (d^\frac12 +1)\bigg[\|b\|_\infty + k\Big(\left\|\sigma\right\|_{\infty} \left\|y\right\|_\beta+ \left[\left\|x_{\cdot-r}\otimes y\right\|_{2\beta (s,t)}+\left\|x_{\cdot-r}\right\|_{\beta(s,t)}\left\|y\right\|_\beta\right]\\
&&\times\left[\left\| \sigma' \right\|_{\infty}+
\left\| \sigma' \right\|_{\gamma} \left\| x_{\cdot-r} \right\|_{\beta(s,t)}^\gamma(t-s)^{\gamma \beta}\right](t-s)^\beta\Big)\bigg].
\end{eqnarray*}
 Now, using the computations that we have done in the first step
\begin{displaymath}\left\| x_{\cdot-r} \right\|_{\beta(s,t)}=\left\| J_1^0 \right\|_{\beta(s-r,t-r)}\le 
\mu (d^\frac12 +1)  \left[1+ 2k+ 3k \left\|y\right\|_{\beta}\right],\end{displaymath}
and (\ref{e5.9}) yields that
\begin{displaymath}\left\|x_{\cdot-r}\right\|_{\beta(s,t)}(t-s)^\beta\le 1.\end{displaymath}
Moreover, also thanks to (\ref{e5.12}) and (\ref{e5.9}) 
\begin{displaymath}\left\|x_{\cdot-r}\otimes y\right\|_{2\beta (s,t)}(t-s)^\beta
=\left\|J^1_2\right\|_{2\beta (s,t)}(t-s)^\beta
\le 1.\end{displaymath}
Combining these two last bounds allow us to show that
\begin{equation}\label{e5.13}
\left\|J_1^1\right\|_{\beta (s,t)}\le \mu  (d^\frac12 +1) \left[1+2k +3k\left\|y\right\|_{\beta}\right].
\end{equation}
\vskip 5pt
\noindent
%%%%%%%%%%%%%%%%%%%%%%%%%
{\underline{Step 2.3:}}  We can  write now our hypothesis of induction. Set  $(H_l)$: for any $s,t \in [lr, (l+1)r]$, the following hypothesis are satisfied 
\begin{displaymath}
(H_l)\qquad \qquad \left\{\begin{array}{l}
\displaystyle \left\|J_1^l\right\|_{\beta (s,t)}\le \mu  (d^\frac12 +1) \left[1+2k+3k \left\|y\right\|_{\beta}\right],\\[2mm]
\displaystyle \left\|J_2^l\right\|_{2\beta (s,t)}\le \mu  (d^\frac12 +1) (1+3k)
\left[ \left\|y\right\|_{\beta}+ \left\|y\right\|_{\beta}^2+
 \left\|y_{\cdot-r}\otimes y\right\|_{2\beta}\right].
\end{array}\right.
\end{displaymath}
We have checked  ($H_1$). Let us suppose that ($H_l$) is satisfied for any $l=1,\dots, n$. Our goal is to proved that $(H_{n+1})$ is also satisfied.
For $s,t \in [(n+1)r,(n+2)r],$, the first two inequalties can be proved as in the previous case:
\begin{eqnarray*}
\left\|J_1^{n+1}\right\|_{\beta (s,t)}&\le &  (d^\frac12 +1) \mu \bigg[1 + k\Big( \left\|y\right\|_\beta+ \left[\left\|x_{\cdot-r}\otimes y\right\|_{2\beta (s,t)}+\left\|x_{\cdot-r}\right\|_{\beta(s,t)}\left\|y\right\|_\beta\right]\\
&&\times\left[1+
 \left\| x_{\cdot-r} \right\|_{\beta(s,t)}^\gamma(t-s)^{\gamma \beta}\right](t-s)^\beta\Big)\bigg],\\
\left\|J_2^{n+1}\right\|_{2\beta (s,t)}&=&\left\|x_{\cdot-r}\otimes y\right\|_{2\beta (s,t)}\\
&\le & \left\|y\right\|_\beta \left\|z\right\|_{\beta(s-r,t-r)}
+\mu \bigg[ \left\|y\right\|_\beta  +k\Big(\left\|y_{\cdot-r}\otimes y\right\|_{2\beta}
+\left\|y\right\|_{\beta}^2\\&&+\left[1+ \left\| x_{.-2r} \right\|_{\beta(s,t)}^\gamma(t-s)^{\gamma \beta}\right]\\
&&\times 
 \Big[\left\| x_{.-2r} \right\|_{\beta(s,t)}\left\|y\right\|_{\beta}^2
 +\left\|x_{.-2r}\otimes y_{\cdot-r}\right\|_{2\beta (s,t)}\left\|y\right\|_{\beta}\\
 &&+\left\| x_{.-2r} \right\|_{\beta(s,t)}
\left\|y_{\cdot-r}\otimes y\right\|_{2\beta}\Big](t-s)^\beta\Big)\bigg].
\end{eqnarray*}
In order to study $J_1^{n+1}$ and $J_2^{n+1}$,  we will use that
\begin{eqnarray*}
\left\|z\right\|_{\beta(s-r,t-r)}&\le & d^\frac12 \mu\left[1+2k+3k\left\|y\right\|_{\beta}\right],\\
\left\| x_{.-2r} \right\|_{\beta(s,t)}(t-s)^\beta&=& \left\|J_1^{n-1}\right\|_{\beta(s-2r,t-2r)}(t-s)^\beta\le 1,\\
\left\|x_{.-2r}\otimes y_{\cdot-r}\right\|_{2\beta (s,t)}(t-s)^\beta&=& 
\left\|x_{\cdot-r}\otimes y\right\|_{2\beta (s-r,t-r)}(t-s)^\beta\\
&=&
\left\|J_2^{n}\right\|_{2\beta(s-r,t-r)}(t-s)^\beta\le 1,\\
\left\|x_{\cdot-r}\right\|_{\beta(s,t)}(t-s)^\beta&=& \left\|J_1^{n}\right\|_{\beta(s-r,t-r)}(t-s)^\beta\le 1.
\end{eqnarray*}
Then, it is not difficult to check that
\begin{eqnarray*}
\left\|J_2^{n+1}\right\|_{2\beta (s,t)}& \le &\mu (d^\frac12 +1) 
(1+3k)
\left[ \left\|y\right\|_{\beta}+  \left\|y\right\|_{\beta}^2+
 \left\|y\otimes y\right\|_{2\beta}\right],\\
\left\|J_1^{n+1}\right\|_{\beta (s,t)}&\le &\mu (d^\frac12 +1) \left[1+2k+3k \left\|y\right\|_{\beta}\right],
\end{eqnarray*}
where the bound of $J_2^{n+1}$ is used to prove the last inequality.
\vskip 5pt
\noindent
%%%%%%%%%%%%%%%%%%%%%%%%%
\hfill $\Box$

%%%%%%%%%%%%%%%%%%%%%%%%%%%%%%%%%%%%%%%%%%%%%%%%%%%%%%%%%%%%%%%%%%%%%%%%%%%
\renewcommand{\theequation}{5.\arabic{equation}} \setcounter{equation}{0}

\section{Stochastic case}

Fix a parameter $H \in(\frac13, \frac{1}{2})$. Set $W^H=\left\{W^H(t),t\in \lbrack 0,T]\right\} $ a $m$-dimensional fractional Brownian
motion of Hurst parameter $H$. The components $W^{H,1},\ldots,W^{H,m}$ are independent  centered Gaussian processes
 with the covariance
function
$$ \ R\left( t,s\right)
=\frac{1}{2}\left( s^{2H}+t^{2H}-|t-s|^{2H}\right).
$$

Let us consider now a Stratonovich type integral with respect to
$W^H$.
Following the approach by Russo and Vallois \cite{RV}, we have:

\begin{definition}
\label{def1}Let $u=\{u(t),t\in \lbrack 0,T]\}$ be a stochastic
process with integrable trajectories. The Stratonovich integral of
$u$ with respect to $W^{H,i}$ is defined as the limit in probability as
$\varepsilon $ tends to zero of
\begin{equation*}
(2\varepsilon )^{-1}\int_{0}^{T}u(s)(W^{H,i}_{s+\varepsilon}
-W^{H,i}_{s-\varepsilon})ds,
\end{equation*}
provided this limit exists. When the limit exists, it is denoted by
$\int_{0}^{T}u(t) \circ dW^{H,i}_t$.
\end{definition}
\vskip 5pt
\noindent
Then, we can consider the tensor product defined by
$$
( W^{H,i}_{\cdot-r} \otimes W^{H,j} )_{s,t} = \int_s^t (W^{H,i}_{v-r} - W^{H,i}_{s-r} ) \circ dW_v^{H,j},$$
for $ 0 \le s < t \le T$ and $i,j=1,\ldots,m.$ Notice that when $i \not= j$, the Stratotovich intergal coincides with the Skorohod integral.

It is proved in Proposition 5.2 in \cite{NNT}  that
\begin{equation}\label{fin1}
E (\vert ( W^{H,i}_{\cdot-r} \otimes W^{H,j} )_{s,t} \vert^p ) \le C_p \vert t- s\vert^{2pH},
\end{equation}
for all $i,j \in \{1,\ldots,m\}$ and $p \ge 1$. Furthermore, they also have checked that there is a version of $( W^{H}_{\cdot-r} \otimes W^{H} )$ such that, for almost
all sample paths of $W^H$, satisfy that $ ( W^{H,i}_{\cdot-r} \otimes W^{H,j} )_{s,t} \in \cc_2^{2\beta}  (\R^{m\times m})$ for any  $\beta \in (\frac13, H)$ and $i,j \in \{1,\ldots,m\},$  where $\cc_2^{2\beta} $ denotes a space of $2 \beta$-H\"older continuous functions of two variables.

Since the multiplicative property can be checked easily using the definition and the properties of the Stratonovich integral, it follows  that $(W^H_{\cdot-r}, W^H, W^H_{\cdot-r}\otimes W^H)$ is a $\beta$-H\"older continuous multiplicative functional for a fixed $\beta \in (\frac13,H)$.

Notice also that (\ref{fin1}) implies that
\begin{equation*}
E (\Vert  W^{H,i}_{\cdot-r} \otimes W^{H,j} )\Vert_{2\beta(r,T)}^p )\le C_p \vert t- s\vert^{2p(H-\beta)},
\end{equation*}
for any $\beta \in (\frac13,H)$ and for all $p \ge 1$.

Then Theorem \ref{tt3} follows easily from  Theorems \ref{tt1}  and Theorem \ref{tt2}, when we apply its results pathwise.

%%%%%%%%%%%%%%%%%%%%%%%%%%%%%%%%%%%%%%%%%%%%%%%%%%%%

\renewcommand{\theequation}{5.\arabic{equation}}
\setcounter{equation}{0}
\section{Appendix}

Let us recall a fixed point Theorem from \cite{B-R}.

\begin{lemma} \label{puntfix}
Let $(X,\rho)$ be a complete metric space, and $\rho_0$ and $\rho_1$ two metrics on X equivalents to $\rho$. If $\cl: X\rightarrow X$ satisfies:
\begin{enumerate}
\item There exists $r_0>0$, $x_0\in X$ such that if $B_0=\left\{x\in X; \rho_0(x_0,x)\leq r_0\right\}$ then $\cl(B_0)\subseteq B_0$,
\item There exists $a\in (0,1)$ such that $\rho_1\left(\cl(x),\cl(y)\right)\leq a\rho_1(x,y)$ for all $x,y\in B_0$.
\end{enumerate}
Then there exists $x^*\in\cl(B_0)\subseteq X$ such that $x^*=\cl(x^*)$.
\end{lemma}

We also need   a result with some properties of the solution of Skorohod's problem.

\begin{lemma} \label{le2}
For each path $\xi\in \cc(\R_+,\R^d)$, there exists a unique solution $(x,z)$ to the Skorokhod problem for $\xi$. Thus there exists a pair of functions  ${\displaystyle (\phi,\varphi): \cc_+(\R_+,\R^d)\rightarrow\cc_+(\R_+,\R^{2d})}$ defined by $\left(\phi(\xi),\varphi(\xi)\right)=(x,z)$. The pair $\left(\phi,\varphi\right)$ satisfies the following:
\vskip 5pt
There exists a constant $K_l>0$ such that for any $\xi_1,\xi_2\in\cc_+(\R_+,\R^d)$ we have for each $t\geq 0$,
\begin{eqnarray*}
\left\|\phi(\xi_1)-\phi(\xi_2)\right\|_{\infty(0,t)}&\leq& K_l\left\|\xi_1-\xi_2\right\|_{\infty(0,t)},\\
\left\|\varphi(\xi_1)-\varphi(\xi_2)\right\|_{\infty(0,t)}&\leq& K_l\left\|\xi_1-\xi_2\right\|_{\infty(0,t)}.
\end{eqnarray*}
Moreover  for each $0\leq s<t\leq T$
\[\left\|\varphi(\xi)\right\|_{\beta(s,t)}\leq d^\frac12 \left\|\xi\right\|_{\beta(s,t)},\]
\end{lemma}
\vskip 6pt
\noindent
{\bf Proof:} We refer the reader to Proposition A.0.1 in \cite{K-W} for the proof of the first part of the Lemma.
\vskip 5pt
\noindent
Set $z=\varphi(\xi)$. Take $u,v$ such that $ s \le u <v \le t$. Fixed a component $i$, we wish to study
$$ \frac{\vert z^i(v) - z^i(u) \vert}{(v-u)^\beta}.$$
When $z^i(v)=z^i(u)$, this is clearly zero.
On the other hand, when $z^i(v)>z^i(u)$, let us define
\begin{eqnarray*}
u^*&:=&\sup\{ u'\ge u; z^i(u)=z^i(u') \},\\
v^*&:=&\inf\{ v' \le v; z^i(v)=z^i(v') \}.
\end{eqnarray*}
Then, $u \le u^* < v^* \le v$ and $z^i(u)=z^i(u^*), z^i(v)=z^i(v^*)$. So
$$ \frac{\vert z^i(v)-z^i(u) \vert}{(v-u)^\beta} \le \frac{\vert  z^i(v^*)-z^i(u^*) \vert}{(v^*-u^*)^\beta} =
\frac{\vert \xi^i(v^*)-\xi^i(u^*) \vert}{(v^*-u^*)^\beta} $$
where the last equality follows from the fact that
$\xi^i$ and $z^i$ coincides whenever $z^i$ is not constant.

Then, note that
$$
\sup_{s < u <v< t} \frac{\vert  z^i(v)-z^i(u) \vert}{(v-u)^\beta}   \le
\sup_{s < u^* <v^*< t} \frac{\vert \xi^i(v^*)-\xi^i(u^*) \vert}{(v^*-u^*)^\beta} \leq
 \Vert \xi \Vert_{\beta(s,t)}.$$
Finally, we get that
$$
\Vert z \Vert_{\beta(s,t)} \leq \Big( \sum_{i=1}^d \big( \sup_{s < u <v < t} \frac{\vert z^i(v)-z^i(u) \vert}{(v-u)^\beta} \big)^2 \Big)^\frac12   \le d^\frac12
 \Vert \xi \Vert_{\beta(s,t)}.$$

\hfill $\Box$

\end{document}